\documentclass[12pt]{article}

\usepackage{hyperref}
\usepackage{xcolor}
\usepackage{amssymb}
\usepackage{amsfonts}
\usepackage{amsthm}
\usepackage{amsmath}
\usepackage{float}
\usepackage{bigints}
\usepackage{fullpage}
\usepackage{mathtools}
\usepackage[utf8]{inputenc} 
\usepackage{cancel}
\usepackage{verbatim}

\newtheorem{theorem}{Theorem}

\newtheorem{lemma}[theorem]{Lemma}

\newtheorem{definition}{Definition}
\newtheorem{problem}[theorem]{Problem}
\newtheorem{corollary}[theorem]{Corollary}
\newtheorem{remark}[theorem]{Remark}

\newcommand{\inter}{\rm int}

\newcommand{\x}{\mathbf{x}}
\newcommand{\p}{\mathbf{p}}
\newcommand{\q}{\mathbf{q}}

\newcommand{\oo}{\mathbf{o}}

\newcommand{\B}{\mathbf{B}}

\newcommand{\A}{\mathbf{A}}

\newcommand{\Ee}{\mathbb{E}}

\newcommand{\Ss}{\mathbb{S}}
\newcommand{\Mm}{\mathbb{M}}

\newcommand{\Ed}{\Ee^d}

\newcommand{\iprod}[2]{\left<#1,#2\right>}

\newcommand{\noshow}[1]{}
\newcommand{\ivol}[2][k]{{\rm V}_{#1}\left(#2\right)}

\title{On a Blaschke-Santal\'o-type inequality for $r$-ball bodies \footnote{\hspace*{.35cm} 2010 Mathematics Subject Classification: 52A20, 52A22.}}

\author{K\'{a}roly Bezdek\thanks{Partially supported by a Natural Sciences and 
Engineering Research Council of Canada Discovery Grant.}
}

\date{}

\begin{document}

\maketitle

\begin{abstract}
Let $\Ee^d$ denote the $d$-dimensional Euclidean space. The $r$-ball body generated by a given set in $\Ee^d$ is the intersection of balls of radius $r$ centered at the points of the given set. The author [Discrete Optimization 44/1 (2022), Paper No. 100539] proved the following Blaschke-Santal\'o-type inequality for $r$-ball bodies: for all $0<k< d$ and for any set of given $d$-dimensional volume in $\Ee^d$ the $k$-th intrinsic volume of the $r$-ball body generated by the set becomes maximal if the set is a ball. In this note we give a new proof showing also the uniqueness of the maximizer. Some applications and related questions are mentioned as well.

\end{abstract}

\section{Introduction}\label{sec:intro}

We denote the Euclidean norm of a vector $\p$ in the $d$-dimensional Euclidean space $\Ee^d$, $d>1$ by $|\p|:=\sqrt{\iprod{\p}{\p}}$, where $\iprod{\cdot}{\cdot}$ is the 
standard inner product. Let $A\subset\Ed$ be a compact convex set, and 
$1\leq k\leq d$. We denote the $k$-th intrinsic volume of $A$ by $\ivol{A}$. It is well known that $\ivol[d]{A}$ is the $d$-dimensional 
volume of $A$, $2\ivol[d-1]{A}$ is the surface area of $A$, and $\frac{2\omega_{d-1}}{d\omega_d}\ivol[1]{A}$ is equal to the mean width of $A$, where $\omega_d$ stands for the volume of a $d$-dimensional unit ball, that is, $\omega_d=\frac{\pi^{\frac{d}{2}}}{\Gamma(1+\frac{d}{2})}$. (For a focused overview on intrinsic volumes see \cite{GHSch}). In this note, for simplicity $\ivol{\emptyset}=0$ for all $1\leq k\leq d$. The closed Euclidean ball of radius $r$ centered at $\p\in\Ed$ is denoted by $\B^d[\p,r]:=\{\q\in\Ed\ |\  |\p-\q|\leq r\}$. 

\begin{definition}\label{r-dual-body}
For a set $\emptyset\neq X\subseteq\Ee^d$, $d>1$ and $r>0$ let the {\rm $r$-ball body} $X^r$ generated by $X$ be defined by $X^r:=\bigcap_{\x\in X}\B^d[\x, r]$. 
\end{definition}

We note that either $X^r=\emptyset$, or $X^r$ is a point, or ${\inter} (X^r)\neq\emptyset$. Perhaps not surprisingly, $r$-ball bodies of $\Ee^d$ have already been investigated in a number of papers however, under various names such as ``\"uberkonvexe Menge'' (\cite{Ma}), ``$r$-convex domain'' (\cite{Fe}), ``spindle convex set'' (\cite{BLNP}, \cite{KMP}), ``ball convex set'' (\cite{LNT}), ``hyperconvex set'' ({\cite{FKV}), and ``$r$-dual set'' (\cite{Be18}). $r$-ball bodies satisfy some basic identities such as $\left((X^r)^r\right))^r=X^r$ and $(X \cup Y)^r=X^r\cap Y^r$, which hold for any $X\subseteq\Ee^d$ and $Y\subseteq\Ee^d$. Clearly, the operation is order-reversing namely, $X\subseteq Y\subseteq\Mm^d$ implies $Y^r\subseteq X^r$. 
In this note we investigate volumetric relations between $X^r$ and $X$ in $\Ee^d$. First, recall the theorem of Gao, Hug, and Schneider \cite{GHSch} stating that for any convex body of given volume in $\Ss^d$ the volume of the spherical polar body is maximal if the convex body is a ball. The author has proved the following Euclidean analogue of their theorem in \cite{Be18}. Let $\A\subset\Ee^d$, $d>1$ be a compact set of volume $V_{d}(\A)>0$ and $r>0$. If $\B\subseteq\Ee^d$ is a ball with $V_{d}(\A)=V_{d}(\B)$, then 
\begin{equation}\label{Bezdek-inequality}
V_{d}(\A^r)\leq V_{d}(\B^r).
\end{equation} 
As the theorem of Gao, Hug, and Schneider \cite{GHSch} is often called a spherical counterpart of the Blaschke--Santal\'o inequality, therefore one may refer to (\ref{Bezdek-inequality}) as a Blaschke--Santal\'o-type inequality for $r$-ball bodies in $\Ee^d$. Next recall that (\ref{Bezdek-inequality}) has been extended by the author to spherical as well as hyperbolic spaces (\cite{Be18}) and then to intrinsic volumes (\cite{Be22}) proving the following Blaschke--Santal\'o-type inequality for intrinsic volumes of $r$-ball bodies in $\Ee^d$ without precise equality condition, which we included here. (See \cite{BeNa} for the core ideas behind Theorem~\ref{Bezdek-inequality-extended} and Theorem 3.1 in \cite{PaPi} for a randomized version.) 

\begin{theorem}\label{Bezdek-inequality-extended}
Let $\A\subset\Ee^d$, $d>1$ be a compact set of volume $V_{d}(\A)>0$ and $r>0$. If $\B\subset\Ee^d$ is a ball with $V_{d}(\A)=V_{d}(\B)$, then 
\begin{equation}\label{Bezdek-inequality-generalized}
V_{k}(\A^r)\leq V_{k}(\B^r)
\end{equation}
holds for all $1\leq k\leq d$ with equality if and only if $\A$ is congruent (i.e., isometric) to $\B$.
\end{theorem}

Fodor, Kurusa, and V\'igh \cite{FKV} have proved the following inequality for $k=d$, which we extended to other intrinsic volumes as well. Corollary~\ref{Fodor-Kurusa-Vigh-extended} follows from Theorem~\ref{Bezdek-inequality-extended}, the homogeneity (of degree $k$) of $k$-th intrinsic volume, and from the observation that $f(x)=x^k(r-x)^k$, $0\leq x\leq r$ has a unique maximum value at $x=\frac{r}{2}$ for any $d>1$, $1\leq k\leq d$, and $r>0$.

\begin{corollary}\label{Fodor-Kurusa-Vigh-extended}
Let $d>1$, $1\leq k\leq d$, $r>0$, and $\A\subset\Ee^d$ be an $r$-ball body. Set $P_k(\A):=V_k(\A)V_k(\A^r)$. Then
\begin{equation}\label{F-K-V}
P_k(\A)\leq P_k\left(\B^d\left[\oo, \frac{r}{2}\right]\right)
\end{equation}
with equality if and only if $\A$ is congruent to $\B^d[\oo, \frac{r}{2}]$.
\end{corollary}
As a further application we mention that Theorem~\ref{Bezdek-inequality-extended} has been used in \cite{Be22} (see also \cite{Be18}) to prove the long-standing Kneser-Poulsen conjecture for uniform contractions of intersections of sufficiently many congruent balls. We close this section with the following complementary question to Theorem~\ref{Bezdek-inequality-extended}, which seems to be new and open, and can be regarded as a Mahler-type problem for $r$-ball bodies.
\begin{problem}\label{Mahler-type}
Let $d>2$, $1\leq k\leq d$, and $0<v<r^d\omega_d=V_d(\B^d[\oo, r])$. Find the minimum of $V_k(\A^r)$ for all $r$-ball bodies $\A\subset\Ee^d$ of given volume $v=V_d(\A)$.
\end{problem}

\begin{remark}\label{2D-Mahler-type}
Problem~\ref{Mahler-type} for $d=2$ can be answered as follows. Let $0<v<\pi r^2$. Then the minimum of $V_1(\A^r)$ (resp., $V_2(\A^r)$) for all $r$-disk domains $\A\subset\Ee^2$ of given area $v=V_2(\A)$ is attained only for $r$-lenses, which are intersections of two disks of radius $r$.
\end{remark}

In the rest of this note we give a short proof for Theorem~\ref{Bezdek-inequality-extended} (which uses the Brunn-Minkowski inequality and the isoperimetric inequality instead of the Alexandrov-Fenchel inequality applied in \cite{Be22}) and derive Remark~\ref{2D-Mahler-type}.

\section{Proof of Theorem~\ref{Bezdek-inequality-extended}}

Clearly, if $\B^r=\emptyset$, then $\A^r=\emptyset$ and (\ref{Bezdek-inequality-generalized}) follows. Similarly, it is easy to see that if $\B^r$ is a point in $\Ee^d$, then (\ref{Bezdek-inequality-generalized}) follows. Hence, we may assume that $\B^r=\B^d[\oo, R]$ and $\B=\B^d[\oo, r-R]$ with $0<R< r$, where $\oo$ denotes the origin in $\Ee^d$.

\begin{definition}
Let $\emptyset\neq K\subset\Ee^d$, $d>1$ and $r>0$. Then the {\rm $r$-ball convex hull} ${\rm conv}_rK$ of $K$ is defined by $${\rm conv}_rK:=\bigcap\{ \B^d[\x, r]\ |\ K\subseteq \B^d[\x, r]\}.$$
Moreover, let the $r$-ball convex hull of $\Ee^d$ be $\Ee^d$. Furthermore, we say that  $K\subseteq\Ee^d$ is {\rm $r$-ball convex} if $K={\rm conv}_rK$.
\end{definition}

\begin{remark}\label{empty}
We note that clearly, ${\rm conv}_rK=\emptyset$ if and only if $K^r=\emptyset$. Moreover, $\emptyset\neq K\subset\Ee^d$ is $r$-ball convex if and only if $K$ is an $r$-ball body.
\end{remark}

We need Lemma 5 of \cite{Be18} stated as

\begin{lemma}\label{basic2}
If $\emptyset\neq K\subseteq\Ee^d$, $d>1$ and $r>0$, then $K^r= ({\rm conv}_rK)^r$.
\end{lemma}
 
Now, via Lemma~\ref{basic2} we may assume that $\A\subset\Ee^d$ is an $r$-ball body of volume $V_{d}(\A)>0$ and $\B=\B^d[\oo, r-R]$ with $0<R< r$ such that $V_{d}(\A)=V_{d}(\B)$. 
Next, recall  Proposition 2.5 of \cite{BeNa} which we state as

\begin{lemma}\label{Bezdek-Naszodi} 
Let $d>1$ and $r>0$. If $\A\subset\Ee^d$ is an $r$-ball body, then $\A+(-\A^r)=\B^d[\oo, r]$, where $+$ denotes the Minkowski sum.
\end{lemma}

Thus, the Brunn-Minkowski inequality for intrinsic volumes (\cite{Gar}, Eq. (74)) and Lemma~\ref{Bezdek-Naszodi}  imply
\begin{equation}\label{B-M-I}
 V_k(\A)^{\frac{1}{k}}+V_k(\A^r)^{\frac{1}{k}}=V_k(\A)^{\frac{1}{k}}+V_k(-\A^r)^{\frac{1}{k}}\leq V_k(\A+(-\A^r))^{\frac{1}{k}}= V_k(\B^d[\oo, r])^{\frac{1}{k}}
\end{equation}
with equality if and only if ($\A$ and $-\A^r$ are homothetic, i.e.,) $\A$ is congruent to $\B$, where $1\leq k\leq d$. Finally, (\ref{B-M-I}), the isoperimetric inequality for intrinsic volumes stating that among convex bodies of given volume the balls have the smallest $k$-th intrinsic volume (\cite{Sc}, Section 7.4), and the homogeneity (of degree $k$) of $k$-th intrinsic volume imply
$$
V_k(\A^r)^{\frac{1}{k}}\leq V_k(\B^d[\oo, r])^{\frac{1}{k}}-V_k(\A)^{\frac{1}{k}}\leq V_k(\B^d[\oo, r])^{\frac{1}{k}}-V_k(\B^d[\oo, r-R])^{\frac{1}{k}}
$$
$$
=V_k(\B^d[\oo, R])^{\frac{1}{k}}=V_k(\B^r)^{\frac{1}{k}},
$$
with $V_k(\A^r)^{\frac{1}{k}}=V_k(\B^r)^{\frac{1}{k}}$ if and only if $\A$ is congruent to $\B$, where $1\leq k\leq d$. This completes the proof of Theorem~\ref{Bezdek-inequality-extended}.

\section{Proof of Remark~\ref{2D-Mahler-type}}
 Let $0<v<\pi r^2$. Let $\A\subset\Ee^2$ be an $r$-disk domain of area $v=V_2(\A)$. Lemma~\ref{Bezdek-Naszodi} implies that
\begin{equation}\label{perimeter}
V_1(\A)+V_1(\A^r)=\pi r.
\end{equation}
Next, according to the reverse isoperimetric inequality of $r$-disk domains, which has been proved by Borisenko and Drach \cite{BoDr} (see \cite{FKV} for a different proof without uniqueness of the extremal set), we have that 
\begin{equation}\label{Bo-Dr}
V_1(\A)\leq V_1(\mathbf{L}),
\end{equation}
where $\mathbf{L}$ denotes the $r$-lens of area $v$ in $\Ee^2$. Moreover, equality holds in (\ref{Bo-Dr}) if and only if $\A$ is conguent to $\mathbf{L}$. Thus, (\ref{perimeter}) and (\ref{Bo-Dr}) imply that 
$\pi r-V_1(\mathbf{L})\leq V_1(\A^r)$. It follows via (\ref{perimeter}) that 
\begin{equation}\label{result-1}
V_1(\mathbf{L}^r)=\pi r-V_1(\mathbf{L})\leq V_1(\A^r)
\end{equation} 
with equality if and only if $\A$ is congruent to $\mathbf{L}$. Finally, observe that (\ref{Bo-Dr}) is equivalent to the statement that if $\A'$ is an $r$-disk domain and $\mathbf{L}'$ is an $r$-lens with $V_1(\A')=V_1(\mathbf{L}')$, then
\begin{equation}\label{dual-Bo-Dr}
V_2(\mathbf{L}')\leq V_2(\A')
\end{equation}
with equality if and only if $\A'$ is conguent to $\mathbf{L}'$. Hence, (\ref{result-1}) combined with (\ref{dual-Bo-Dr}) yields
\begin{equation}\label{result-2}
V_2(\mathbf{L}^r)\leq V_2(\A^r)
\end{equation} 
with equality if and only if $\A$ is congruent to $\mathbf{L}$. This completes the proof of Remark~\ref{2D-Mahler-type}.

\bigskip


\noindent K\'aroly Bezdek \\
\small{Department of Mathematics and Statistics, University of Calgary, Canada}\\
\small{Department of Mathematics, University of Pannonia, Veszpr\'em, Hungary\\
\small{E-mail: \texttt{bezdek@math.ucalgary.ca}}


\small
\begin{thebibliography}{GGM}







\bibitem{BLNP}
K. Bezdek, Zs. L\'angi, M. Nasz\'odi, and P. Papez, Ball-polyhedra, {\it Discrete Comput. Geom.} {\bf 38/2} (2007), 201--230.



\bibitem{Be18}  K. Bezdek, From r-dual sets to uniform contractions, {\it Aequationes Math.} \textbf{92/1} (2018), 123--134.

\bibitem{BeNa}
K. Bezdek and M. Nasz\'odi, The Kneser-Poulsen conjecture for special contractions, {\it Discrete Comput. Geom.} \textbf{60/4} (2018), 967--980.

\bibitem{Be22}
K. Bezdek, On the intrinsic volumes of intersections of congruent balls, {\it Discrete Optim.} \textbf{44/1} (2022), Paper No. 100539 (7 pages).

\bibitem{BoDr}
A.~A. Borisenko and K.~D.~Drach, Isoperimetric inequality for curves with curvature bounded below, {\it Translation of Mat. Zametki} \textbf{95/5} (2014), 656--665, {\it Math. Notes} \textbf{95/5-6} (2014), 590--598.

  






\bibitem{Fe}
L. Fejes T\'oth, Packing of r-convex discs, {\it Studia Sci. Math. Hungar.} {\bf 17/1-4} (1982), 449--452.

\bibitem{FKV}
F. Fodor, \'A. Kurusa, and V. V\'igh, Inequalities for hyperconvex sets, {\it Adv. Geom.} {\bf 16/3} (2016), 337--348.

\bibitem{GHSch}
F. Gao, D. Hug, and R. Schneider, Intrinsic volumes and polar sets in spherical space, {\it Math. Notae} {\bf 41} (2003), 159--176. 

\bibitem{Gar}
R. J. Gardner, The Brunn-Minkowski inequality, {\it Bull. Am. Math. Soc.} {\bf 39/3} (2002), 355--405.









\bibitem{KMP}
Y. S. Kupitz, H. Martini, and M. A. Perles, Ball polytopes and the V\'azsonyi problem, {\it Acta Math. Hungar.} {\bf 126/1-2} (2010), 99--163.


\bibitem{LNT}
Zs. L\'angi, M. Nasz\'odi, and I. Talata, Ball and spindle convexity with respect to a convex body, {\it Aequationes Math.} {\bf 85/1-2} (2013), 41--67.

\bibitem{Ma}
A. E. Mayer, Eine \"Uberkonvexit\"at, {\it Math. Z.} {\bf 39/1} (1935), 511--531.

\bibitem{PaPi}
G. Paouris and P. Pivovarov, Random ball-polyhedra and inequalities for intrinsic volumes, {\it Monatsh. Math.} {\bf 182/3} (2017), 709--729. 


\bibitem{Sc} R. Schneider, {\it Convex bodies: the Brunn-Minkowski theory}, Encyclopedia of Mathematics and its Applications, vol. 44, Cambridge University Press, Cambridge, 1993.






\end{thebibliography}
\end{document}